\documentclass[12pt]{article}
\usepackage{fullpage,amsfonts,amsmath,amsthm,graphicx}

\usepackage[usenames]{color}

\begin{document}

\newcommand{\mm}[1]{{\color{black}{#1}}}

\newcommand{\mmm}[1]{{\color{black}{#1}}}

\def\a{\alpha}
\def\b{\beta}
\def\c{\chi}
\def\d{\delta}
\def\D{\Delta}
\def\e{\varepsilon}
\def\f{\phi}
\def\F{\Phi}
\def\g{\gamma}
\def\G{\Gamma}
\def\K{\Kappa}
\def\z{\zeta}
\def\th{\theta}
\def\Th{\Theta}
\def\la{\lambda}
\def\La{\Lambda}
\def\m{\mu}
\def\n{\nu}
\def\p{\pi}
\def\P{\Pi}
\def\r{\rho}
\def\R{\Rho}
\def\s{\sigma}
\def\S{\Sigma}
\def\t{\tau}
\def\om{\omega}
\def\Om{\Omega}
\def\smallo{{\rm o}}
\def\bigo{{\rm O}}
\def\to{\rightarrow}
\def\E{{\bf E}}
\def\ex{{\bf E}}
\def\cd{{\cal D}}
\def\rme{{\rm e}}
\def\hf{{1\over2}}
\def\R{{\bf  R}}
\def\cala{{\cal A}}
\def\cale{{\cal E}}
\def\Fscr{{\cal F}}
\def\cc{{\cal C}}
\def\calc{{\cal C}}
\def\calh{{\cal H}}
\def\call{{\cal L}}
\def\calr{{\cal R}}
\def\calb{{\cal B}}
\def\calz{{\cal Z}}
\def\calq{{\cal Q}}
\def\bk{\backslash}

\def\out{{\rm Out}}
\def\temp{{\rm Temp}}
\def\keep{{\rm Keep }}
\def\overused{{\rm Overused}}
\def\big{{\rm Big}}
\def\notbig{{\rm Notbig}}
\def\moderate{{\rm Moderate}}
\def\swappable{{\rm Swappable}}
\def\candidate{{\rm Candidate}}
\def\bad{{\rm Bad}}
\def\crit{{\rm Crit}}
\def\col{{\rm Col}}
\def\dist{{\rm dist}}
\def\eq{{\rm Eq}}
\def\uq{{\rm Uq}}
\def\poly{{\rm poly}}

\newcommand\fix{{\rm FIX }}
\newcommand\fixx{{\rm FIX2 }}
\newcommand{\blank}{{\mathsf{Blank}}}

\newcommand{\Exp}{\mbox{\bf E}}
\newcommand{\Med}{\mbox{\bf Med}}
\newcommand{\med}{\mbox{\bf Med}}
\newcommand{\var}{\mbox{\bf Var}}
\newcommand{\pr}{\mbox{\bf Pr}}
\newcommand{\xmax}{X_{\rm max}}

\newtheorem{lemma}{Lemma}
\newtheorem{theorem}[lemma]{Theorem}
\newtheorem{corollary}[lemma]{Corollary}
\newtheorem{claim}[lemma]{Claim}
\newtheorem{remark}[lemma]{Remark}
\newtheorem{observation}[lemma]{Observation}
\newtheorem{proposition}[lemma]{Proposition}
\newtheorem{definition}[lemma]{Definition}

\newcommand{\limninf}{\lim_{n \rightarrow \infty}}
\newcommand{\proofstart}{{\bf Proof\hspace{2em}}}
\newcommand{\tset}{\mbox{$\cal T$}}
\newcommand{\proofend}{\hspace*{\fill}\mbox{$\Box$}\vspace{2ex}}
\newcommand{\bfm}[1]{\mbox{\boldmath $#1$}}
\newcommand{\reals}{\mbox{\bfm{R}}}
\newcommand{\expect}{\mbox{\bf Exp}}
\newcommand{\he}{\hat{\e}}
\newcommand{\card}[1]{\mbox{$|#1|$}}
\newcommand{\rup}[1]{\mbox{$\lceil{ #1}\rceil$}}
\newcommand{\rdn}[1]{\mbox{$\lfloor{ #1}\rfloor$}}
\newcommand{\ov}[1]{\mbox{$\overline{ #1}$}}
\newcommand{\inv}[1]{\frac{1}{ #1}}

\def\calc{{\cal C}}
\def\cald{{\cal D}}

\title{Asymptotically good edge correspondence colourings}
 
\author{Michael Molloy\thanks{Dept of Computer Science,
University of Toronto, molloy@cs.toronto.edu.  Research supported by  NSERC Discovery Grant 2019-06522.}
 \and Luke Postle\thanks{Dept of Combinatorics and Optimization,
University of Waterloo, lpostle@uwaterloo.ca. Research supported by NSERC Discovery Grant 2019-04304 and the Canada Research Chairs program.}
}

\maketitle

\begin{abstract}  We prove that \mm{every} simple graph with maximum degree $\D$ has  edge correspondence number $\D+o(\D)$.
\end{abstract}

\section{Introduction}

Graph colouring is one of the richest and most fundamental fields of graph theory.  In its most basic form, one must assign a colour from a given set to each vertex of a graph \mm{so} that the endpoints of each edge \mm{get} different colours.  Many variations have arisen, one of the most fruitful being list colouring: \mm{Vizing \cite{vlist} and  Erd\H{o}s, Rubin Taylor~\cite{ert}} independently suggested that rather than assigning colours to all vertices from a single set, we can give each vertex $v$  its own list of permissable colours, $L(v)$.  This very natural variation grew into a prominent subfield of graph colouring.

Recently, Dvo\u{r}\'{a}k and Postle~\cite{dp} introduced another natural variation, {\em correspondence colouring}. Rather than using the same colouring rule for all edges, each edge \mm{can forbid a different set of pairs of colours on its endpoints. The only requirement is that no colour can be forbidden to a vertex by two pairs on the same edge}. Specifically, each edge $uv$ is given a partial matching $M_{uv}$ between $L(u),L(v)$.  The goal is to assign to each vertex $v$ a colour from $L(v)$ \mm{so} that for every edge $uv$, the colours assigned to $u$ and $v$ are not paired in $M_{uv}$.  Note that if every edge $uv$ simply matches each colour in $L(u)\cap L(v)$ to itself then we have the usual list colouring.  \mm{Several studies of correspondence colouring have already appeared}, eg.\ \cite{ab1,ab2,bk1,bk,akz,fh,ko}.  \mm{See~\cite{bk} for a discussion of how correspondence colouring can be more challenging than list colouring, including how some common useful approaches to list colouring do not apply to correspondence colouring.}

\mm{In an instance of correspondence colouring, if every list of colours has the same size, $k$, then we can assume that each list is $\{1,...,k\}$.  To see this, consider an instance where the lists differ. For every vertex $v$ take a bijection $\s_v:L(v)\rightarrow\{1,...,k\}$ and for every edge $uv$, replace each $(i,j)\in M_{uv}$ with $(\s_v(i),\s_v(j))$. Similarly, when the lists have different sizes, we can assume that the list of each vertex $v$ is $\{1,...,|L(v)|\}$. So there is no difference between e.g.\ the correspondence number and the list correspondence number of a graph.}

One of the most pursued open questions in list colouring is: \mm{When edge-colouring a simple graph (i.e. assigning colours to the edges so that every two edges which share a vertex must get different colours), are the identical lists the most difficult lists?}  In other words, is the list edge chromatic number of a simple graph equal to the edge chromatic number?   This has been answered  in the affirmative for specific classes of graphs  (eg.\ \cite{fg,pw}), but is still open for general graphs.  In a seminal paper~\cite{jk}, Kahn proved that the two numbers are asymptotically equal: the list edge chromatic number of a simple graph with maximum degree $\D$ is equal to $\D+o(\D)$. In a followup paper~\cite{jk2} he proved that the two numbers are asymptotically equal for multigraphs as well.  Molloy and Reed~\cite{mrlcc} showed that, for simple graphs, the $o(\D)$ term is at most $\sqrt{\D}\poly(\log\D)$.   See~\cite{jt} for a more thorough background to list colouring.

\mm{Correspondence colouring can be defined for edge colouring in a natural way: each pair of edges that share a vertex is given a list of forbidden pairs (this is defined more formally below).}
Bernshteyn and Kostochka~\cite{bk} showed that  the edge correspondence  number of a simple graph can exceed the edge chromatic number.   In fact, every $\D$-regular simple graph has  edge correspondence number at least $\D+1$, whereas many such graphs have edge chromatic number $\D$.   However,  we show here that Kahn's result holds in this context; i.e.\ every simple graph with maximum degree $\D$ has edge correspondence  number  $\D+o(\D)$.  \mm{The previous best bound in this direction was  $(2-\e)\D$ for a constant $\e>0$, which follows from work in~\cite{bpp} (in particular, the correspondence colouring version of their Theorem 1.6).}

To set things up formally: We are given a simple graph $G$ and a set of colours $\calq=\{1,...,q\}$.  For each pair of incident edges $e,f$, we are given a partial matching $M_{e,f}$ on $(\calq,\calq)$; i.e. a collection of at most $q$ pairs $(\a,\a')\in\calq\times\calq$ such that each colour $\a$ is the \mm{first} element of at most one pair and the \mm{second} element of at most one pair. $M_{f,e}$ will consist of the reversal of all pairs in $M_{e,f}$, so there is only one matching on each pair of incident edges. This collection of partial matchings is called an {\em edge correspondence}.  An {\em edge correspondence colouring} is an assignment to each edge $e\in E(G)$ of a colour  $\s(e)\in\calq$ , such that for every two incident edges $e,f$, the pair $(\s(e),\s(f))$ \mm{is not in} $M_{e,f}$. 

The {\em edge correspondence number} of a graph $G$ is the minimum $q$ such that an edge correspondence colouring exists for every edge correspondence.  We denote this \mm{by} $\chi'_{DP}(G)$, following the notation of Bernshteyn and Kostochka who refer to correspondence colouring as {\em DP-colouring}, using the intials of the founders.

\begin{theorem}\label{mt} Let $G$ be any simple graph with maximum degree $\D$.  Then $\chi'_{DP}(G)=\D+o(\D)$.
\end{theorem}

\mm{{\bf Remark:} Throughout the paper, asymptotic notation is with respect to $\D\rightarrow\infty$.

We colour the edges using an iterative procedure, first introduced in Kahn's proof~\cite{jk} and since then adapted to a very large number of results (see eg.~\cite{mrbook}).  At each step, we colour a small proportion (roughly $\inv{\ln\D}$) of the edges.  We do so by considering a random colour assignment to those edges.  If we name a particular vertex, then a probabilistic analysis shows that the colours of the edges near that vertex will likely satisfy certain properties, for example that each edge has many remaining colours that can still be legally assigned to it.  We  apply the Lovasz Local Lemma to obtain a colouring in which the colours of the edges near {\em every} vertex satisfies those properties.  Eventually we will have coloured almost all the edges; the remaining edges will be such that they are easily dealt with.

One useful aspect to our procedure: we carry out the random colouring so that for any edge $e=uv$, the effect of random choices involving edges incident to $u$ is independent of the effect of choices on edges incident to $v$, and we track the cummulative affects of those choices seperately.  This is where we make critical use of the facts: (i) $G$ is simple, and (ii) edge colouring has the nice structural property that the neighbourhood of each edge consists of two disjoint cliques.
}

\section{Preliminaries} 
\subsection{Setup}\label{s2.1}

Let $\e$ be any sufficiently small constant. We will prove that there exists $\D(\e)$ such that if $\D(G)\geq \D(\e)$ then $\chi_{DP}'(G)\leq(1+\e)\D(G)$.  This is enough to establish Theorem~\ref{mt}. We do not name $\D(\e)$ explicitly; instead we just assume that $\D(G)$ is large enough to satisfy various inequalities that depend on $\e$.

So we are given a graph $G$ with maximum degree $\D$, colours $\calq=\{1,...,(1+\e)\D\}$ and an edge correspondence.  Our goal is to prove that, so long as $\D$ is sufficiently large in terms of $\e$, there must be an edge correspondence colouring.

\mmm{
\begin{definition} For two incident edges $e,f$ and colours $\a,\a'$, we say that $\a:e$ {\em blocks} $\a':f$ if $(\a,\a')$ is an edge of $M_{e,f}$; i.e. if we are not permitted to assign $\a$ to $e$ and assign $\a'$ to $f$.
\end{definition}

Note that  $\a:e$  blocks $\a':f$ iff  $\a':f$  blocks $\a:e$.
}

\subsection{Probabilistic tools}\label{spt}

 We often use the following straightforward bound:
\[{a\choose b}\leq \left(\frac{ea}{b}\right)^b.\]

We also rely on the following standard tool of the probabilistic method.  

{\bf The Lov\'asz Local Lemma} \cite{el1}. {\em Let $\cala=\{A_1,...,A_n\}$
be a set of random events \mm{such} that for each $1\leq i\leq n$:
\begin{enumerate}
\item[(i)] $\pr(A_i)\leq p$; and
\item[(ii)] $A_i$ is mutually independent
of all but at most $d$ other events.
\end{enumerate}
If $pd\leq\inv{4}$ then $\pr(\ov{A_1}\cap...\cap\ov{A_n})>0$.
}

$BIN(n,p)$ is the sum of $n$ independent random \mm{boolean} variables
where each is equal to 1 with probability $p$. The following is a simplified special case of Chernoff's 
bound.  It follows from, e.g. Corollary A.1.10 and Theorem A.1.13 from
Appendix A of \cite{as}.

{\bf The Chernoff Bound.}  {\em For any $0<t\leq np$:
\[\pr(|BIN(n,p)-np|>t)<2e^{-t^2/3np}.\]
}

Theorem 2.3 from~\cite{cmsurv} generalizes the Chernoff Bound. Parts (b,c) of that theorem imply: 

\begin{lemma}\label{l23}
Suppose that we have independent random variables $Z_1,...,Z_n$, with $0\leq Z_i\leq 1$ for each $i$.  Set $Z=\sum_{i=1}^n Z_i$. For any $0<t\leq \ex(Z)$:

\[\pr(|Z-\ex(Z)|>t)<2e^{-t^2/3\ex(Z)}.\]
\end{lemma}

Our final concentration tool is Talagrand's Inequality, which often provides a stronger bound when \mm{the expectation of a random variable is
much smaller than the number of trials that determine it}.  We will use the following variant on Talagrand's original statement  from~\cite{mt}.  The proof is deferred to an appendix.

\mmm{
\noindent{\bf Talagrand's Inequality.}
{\em Let $X$ be a random variable determined by the
independent trials $T_1,...,T_n$. Let $\xmax,D>0$ and suppose that $X$ always satisfies $0\leq X\leq\xmax$. Let $F$ be the event that for the outcome $y=(y_1,...,y_n)$ of the trials there exists $b_1,....,b_n>0$ such that
\begin{enumerate}
\item[(i)] $\sum_{i=1}^n b_i^2\leq D$; and
\item[(ii)] for any possible outcome $z=(z_1,...,z_n)$ of the trials, we have
\[X(z)\geq X(y)-\sum_{y_i\neq z_i} b_i.\]
\end{enumerate}
Then for any $t\geq 0$ we have
\begin{equation}\label{tigoal}
\pr\left(|X-\ex(X)|>t +  35\sqrt{D}+\xmax\times\pr(\ov{F})\right)\leq 2\pr(\ov{F}) + 5e^{-\frac{t^2}{4D}}.
\end{equation}
}
}

\subsection{Adapting previous work}  We will find an edge correspondence colouring of the given graph using a common randomized procedure.  One feature of that procedure is that when a edge gets a colour, any conflicting colours are removed from the lists of available colours for all neighbouring edges.

We would like to have applied the argument from~\cite{mrlcc} to prove that $\chi'_{DP}(G)=\D+\sqrt{\D}\poly(\log\D)$.  The hurdle we could not overcome is as follows:  The procedure in \cite{mrlcc} begins by reserving a set of colours \mm{at each vertex which cannot be assigned to any  edges incident to that vertex.   In the context of list edge colouring, this ensures that at the end of the procedure, each uncoloured edge $uv$ can be assigned any of the colours that were reserved at both $u$ and $v$; however, this is not true for correspondence colouring.}

So instead we followed what is, at heart, the argument from~\cite{jk}, although presented as in~\cite{mrlcc,mrbook}.  One difference is as follows:  in that argument, one kept track of a parameter $T(v,c)$ which was the set of edges incident to $v$ which could still receive the colour $c$.  That parameter was important because $T(u,c)$ and $T(v,c)$ comprised the edges which could cause $c$ to be removed from the list of the edge $e=uv$.  In the context of correspondence colouring, we need to redefine that parameter.  For each edge $e=uv$ we define $T(e,v,c)$ to be the set of edges incident to $v$ which can still receive a colour that will cause $c$ to be removed from the list of $e$.  Once that parameter is defined, the remainder of the argument is a simple adaptation of those from~\cite{jk,mrlcc}.

\subsection{A quick result}
The following result is a very simple variation on the main result of~\cite{brlist} (improved in~\cite{ph}).  It will be used at the end of our proof, just as the result of~\cite{brlist} is used at the end of many similar proofs.

As mentioned above, one can assume that each edge has the same list of permissible colours.  Nevertheless, it will be convenient to extend the definition of an edge correspondence, and an edge correspondence colouring in the obvious way to the case where the lists may differ.  

\mm{We use $f\sim e$ to denote that edges $f,e$ are adjacent. For a pair of incident edges $e,f$, we say that $\a\in L(e)$ has a partner in $M_{e,f}$ if $\a$ is the first element of one of the pairs in $M_{e,f}$; i.e. if assigning $\a$ to $e$ forbids a colour to be assigned to $f$.}

\begin{lemma}\label{ldpcc}
We are given  a simple graph $G$;  a list $L(e)$ of size at least $L$ on each edge $e$; and an edge correspondence such that for each edge $e$ and colour $\a\in L(e)$, there are at most $T$ edges $f\sim e$ such that $\a$ has a partner in $M_{e,f}$.  If $L\geq 8T$ then there is an edge correspondence colouring.
\end{lemma}

\mm{The proof is essentially identical to that from~\cite{brlist}; we include it here for completeness.  It also follows easily from Theorem 2 of~\cite{ph}, with the constant 8 improved to 2.}

\proofstart  
Assign to each edge $e$ a uniformly random colour from $L(e)$.  For each pair of incident edges $e,f$ and pair of colours $(\a,\a')\in M_{e,f}$ we define $A_{e,f,\a,\a'}$ to be the event that $e$ is assigned $\a$ and $f$ is assigned $\a'$.  The probability of each \mm{such} event is at most $1/L^2$. \mm{Each event $A_{e,f,\a,\a'}$ is easily seen to be mutually independent of all events which do not involve $e$ or $f$}; i.e. of all but at most $2LT$ other events.  Since $\inv{L^2}\times 2LT\leq \inv{4}$, the Lov\'asz Local Lemma implies that with positive probability none of these events hold; i.e. we obtain an edge correspondence colouring.
\proofend

\section{A random colouring procedure}\label{srcp}
We colour the graph randomly through a series of iterations, as described in the introduction.
Roughly speaking, at each iteration we colour a small proportion of the edges. When an edge receives a colour then we remove any conflicting colours from the lists of incident edges. If two incident edges receive conflicting colours then both are uncoloured.  A few technical clarifications: 

(a) When an edge $e$ receives a colour then conflicting colours are removed \mm{from all incident edges {\em even if that colour is removed from $e$}.}  This is often refered to as {\em wasteful} since some colours are needlessly removed
from lists.  We do this because it simplifies the analysis.   Furthermore, because such a small proportion of edges are coloured, a vanishing proportion of coloured edges \mm{have their colour removed}.  As a result, the number of colours removed needlessly from each list is negligible.

(b) We allow each edge to receive multiple colours. For each edge $e$ and each colour $c\in L(e)$, the current list for $e$, we assign $c$ to $e$ with probability $1/(|L(e)|\ln\D)$; the choice of whether to assign $c$ to $e$ is independent of the choices for all other colours in $L(e)$.  So the probability that $e$ gets at least one colour is roughly $1/\ln \D$.  Making these assignments independently simplifies the analysis.  And the probability that $e$ gets at least two colours is $O(\ln^{-2}\D)$ which is small enough to be negligible.  We believe that this technique was first used by Johansson in~\cite{aj}.

(c) It is very convenient if, at each iteration, all lists have the same size and the probability that a colour $c$ is removed from $L(e)$ is the same for every $c,e$.  We enforce this by truncating some lists and by carrying out so-called {\em equalizing coin flips} which round up the probability of a colour being removed from a list.  

\mm{Our procedure makes use of the parameters $L_i,T_i, \eq_i(e,c)$.  They will be defined formally below, as their definitions will be more intuitive after reading the procedure. For now, the main things to understand are: (i) our analysis will enforce that at the beginning of each iteration $i$, every edge $e$ has $|L(e)|\geq L_i$; (ii) $\eq_i(e,c)$ is the value required for the equalizing coin flips described above. 

If, during a particular iteration, colour $c$ is assigned to $e$  (in Step 2(b.i)) and colour $c$ is not unassigned from $e$ (in Step 2(b.ii) or Step 2(c)) then we say that $e$ {\em retains} $c$.   At any iteration, an edge is considered {\em uncoloured} if it did not retain any colour during the previous iterations.

Recall from Section~\ref{s2.1} that $\calq=\{1,...,(1+\e)\D\}$. }

\begin{enumerate}
\item Initialize for every edge $e=uv$ and colour $c$: $L(e)=\calq$, $T(e,v,c)$ is the set of all edges incident to $e$ at $v$.
\item For each $i\geq 1$ until \mm{$L_i<\D^{9/10},T_i<\D^{9/10}$ or $L_i>10 T_i$}:
\begin{enumerate}
\item For every uncoloured edge $e$ with $|L(e)|>L_i$, remove $|L(e)|-L_{i}$ arbitrary colours from $L(e)$.
\item For every uncoloured edge $e$ and every colour $c\in L(e)$:
\begin{enumerate}
\item assign $c$ to $e$ with probability $1/(L_i\ln \D)$.
\item If $c$ \mm{was} assigned to $e$ then for every $f\sim e$, if there is a colour $c'\in L(f)$ with $c':f$ blocking $c:e$ then
\begin{enumerate}
\item  remove $c'$ from $L(\mm{f})$; and
\item if $c'$ was assigned to $f$ then unassign $c'$ from $f$
\end{enumerate}
\end{enumerate}
\item For every colour $c$ still in $L(e)$, with probability $1-\eq_i(e,c)$: remove $c$ from $L(e)$, and if $c$ was assigned to $e$ then unassign $c$ from $e$. 
\end{enumerate}
\end{enumerate}

When this procedure terminates, each edge that has \mm{is not  uncoloured is given one of the colours that it retained.  We will argue that it will terminate because $L_i>10T_i$ which will imply that} this partial edge correspondence colouring can be completed using Lemma~\ref{ldpcc}.

\mm{Consider an edge $e=uv$. As mentioned in the introduction, we wish to seperate the random choices related to the effect on $e$ of edges around $v$ from the effect of the edges around $f$. So to carry out the choice in line~2(c) whether to keep $c$ in $L(e)$, we will in fact make two independent random coin flips $F(e,u,c),F(e,v,c)$, which return 1 with probabilities $\eq_i(e,u,c)$ and $\eq_i(e,v,c)$, respectively.  If either returns 0 then $c$ is removed from $L(e)$.  (The values of those probabilities are specified below.)

For each edge $e=uv$ and colour $c$, we define the following sets  at the beginning of step 2 of iteration $i$, i.e. after the lists have been truncated:

\begin{eqnarray*}
L_i(e)&=&\mbox{ the set of colours remaining in the list on } e\\ 
T_i(e,v,c)&=&\mbox{ the set of \mm{uncoloured} edges $f$ containing $v$ for which}\\
&&\qquad\mbox{ there is a colour $c'\in L_i(f)$ such that $c': f$ blocks $c: e$.}
\end{eqnarray*}

We will recursively define parameters $L_i,T_i$ and enforce that for each iteration $i$:
\begin{equation}\label{etrunc}
|L_i(e)|= L_i \mbox{ and } |T_i(e,v,c)|\leq T_i \mbox{ for every edge $e$, \mm{ endpoint $v$ of $e$} and colour $c\in L_i(\mm{e})$.}
\end{equation}
Note that the first condition means that $|L(e)|\geq L_i$ at the beginning of iteration $i$.

\mm{Recalling that we wish to focus seperately on colours removed from $L(e)$ because of edges around $u$ and those removed because of colours around $v$, we introduce the following terminology:

\begin{definition} For an edge $e=uv$ with $c\in L(e)$.  We say that {\em $L(e)$ loses $c$ at $v$} during iteration $i$ if either (a) some edge $f$ with endpoint $v$ is assigned a colour $c'$ where $c':f$ blocks $c:e$ or (b) the equalizing coin flip $F(e,v,c)$ returns 0.
\end{definition}

Note that if $e$ is assigned the colour $c$ in step~2(b) then $c$ is unassigned from $e$  iff  $L(e)$ loses $c$ at $u$ or $L(e)$ loses $c$ at $v$.

 Suppose that~(\ref{etrunc}) holds at the beginning of step 2 of iteration $i$.   Thus the probability that  no colour $c'$ is assigned to an edge $f=wv$ where $c':f$ blocks $c:e$, is  $\left(1-\inv{L_i\ln\D}\right)^{|T_i(e,v,c)|}\geq\left(1-\inv{L_i\ln\D}\right)^{2T_i}$.  }

This inspires us to define
\begin{eqnarray*}
\keep_i&=&\left(1-\inv{L_i\ln\D}\right)^{T_i}\\
\eq_i(e,u,c)&=&\keep_i/\left(1-\inv{L_i\ln\D}\right)^{|T(e,u,c)|}\\
\eq_i(e,v,c)&=&\keep_i/\left(1-\inv{L_i\ln\D}\right)^{|T(e,v,c)|}\\
\eq_i(e,c)&=&\eq_i(e,u,c)\times \eq_i(e,v,c)
\end{eqnarray*}

So the probability that $L(e)$ loses $c$ at $v$ during iteration $i$ is exactly $1-\keep_i$, and the event that it loses $c$ at $v$ is independent of the event that it loses $c$ at $u$ (since the graph is simple).  

Thus, the probability that $c$ remains in $L(e)$ at the end of iteration $i$ is exactly $\keep^2_i$, and so
 the expected number of such colours remaining  on $L(e)$ is $L_i\times \keep^2_i$. 
}

We now turn our attention to $T_{i+1}(e,v,c)$. We cannot show that this parameter is concentrated because it is possible for the assignment of a single colour to some $f\sim e$ to cause $T_{i+1}(e,v,c)$ to drop to $\emptyset$.  So instead, we focus on a related parameter which essentially removes the influence of edges incident to $v$.

$T'_{i+1}(e,v,c)$ is defined to be the set of edges $f=vw\in T_{i}(e,v,c)$ such that (a) $f$ does not retain a colour during iteration $i$ and (b) $L(f)$ does not lose $c'$  at $w$ during iteration $i$, where $c'$ is the unique colour in $L(f)$ such that $c':f$ blocks $c:e$.

\mm{Note that $T_{i+1}(e,v,c)\subseteq T'_{i+1}(e,v,c)$.} So an upper bound on $|T'_{i+1}(e,v,c)|$ will provide an upper bound on $|T_{i+1}(e,v,c)|$.  The fact that each colour in the list of an edge is assigned to that edge independently, makes it  simple to bound the expectation of $|T'_{i+1}(e,v,c)|$:

For any edge $f=vw$ and any colour $\a\in L_i(f)$, let $Z(\a,f)$ be the event that $\a$ is assigned to $f$ and let $Y^v(\a,f), Y^w(\a,f)$ be the events that $L(f)$ loses $\a$ at $v$, and $L(f)$ loses $\a$ at $w$ during iteration $i$.  The following observation is very helpful:

\begin{observation}\label{oZY}  The events $\{ Z(\a,f), Y^v(\a,f), Y^w(\a,f): \a\in L_i(f)\}$  are mutually \mm{independent.}
\end{observation}
\proofstart  First, by the way we carry out Step 2(b), the events $Z(c,e)$  over all edges $e$ and $c\in L_i(e)$ are  determined by independent trials. 
$Y^v(\a,f)$ is determined by the events $Z(\mm{h},\a')$ for all edges $\mm{h}\in T(f,v,\a)$ and colours $\a'\in L(g)$ such that $\a':\mm{h}$ blocks $\a:f$.   By the nature of correspondence colouring, $\a':\mm{h}$ can block $\a:f$ for at most one colour $\a$.  Since the graph is simple, no edge $\mm{h}$ is relevant to both a $Y^v(\cdot ,f)$ event and a $Y^w(\cdot ,f)$ event.  So these events are determined by disjoint sets of trials.
\proofend

Now consider any $f\in T_i(e,v,c)$, where $c':f$ blocks $c:e$. \mm{Suppose that~(\ref{etrunc}) holds at the beginning of iteration $i$. Then} Observation~\ref{oZY} implies (see explanation below):
\begin{eqnarray*}
\pr\left(f\in T'_{i+1}(e,v,c)\right)&=&\keep_i\times\left(1-\inv{\ln\D L_i}\keep_i\right)
\times\prod_{\a\in L_i(f), \a\neq c'}\left(1-\inv{\ln\D L_i}\keep^2_i\right)\\
&<&\mm{\left(1-\inv{\ln\D L_i}\keep^2_i\right)^{L_i}\qquad\qquad\qquad\qquad \mbox{since $\keep_i<1$}}\\
&<&\keep_i\times\left(1-\frac{\mm{1-\e/2}}{\ln\D}\keep^2_i\right).
\end{eqnarray*}
{\bf Explanation:} The first term is the probability that $L(f)$ does not lose $c'$ at $w$.  The second term is the probability that if $f$ is assigned $c'$ then $L(f)$ loses $c'$ at $v$ and so $c$ is removed from $f$.  The third term is the probability that $f$ does not retain any other colour.

\mm{This yields that if~(\ref{etrunc}) holds for iteration $i$ then}:
\begin{equation}\label{et'}\ex[ |T'_{i+1}(e,v,c)|]< |T_{i}(e,v,c)| \times\left(1-\frac{\mm{1-\e/2}}{\ln\D}\keep^2_i\right)\times\keep_i.
\end{equation}

We will prove in section~\ref{sconc} that   $|T'_{i+1}(e,v,c)|$ and the \mm{number of} colours removed from $L(e)$ during step 2 are both concentrated.  This leads us to recursively define: $L_0=(1+\e)\D, T_0=\D$ and
\begin{eqnarray}
L_{i+1}&=& L_i\times \keep^2_i - \D^{2/3}\label{erec1}\\
T_{i+1}&=& T_i\times\left(1-\frac{\mm{1-\e/2}}{\ln\D}\keep^2_i\right)\times\keep_i + \D^{2/3}.\label{erec2}
\end{eqnarray}

\mm{{\bf Remark:} Recall that our procedure halts if $L_i$ or $T_i$ drops below $\D^{9/10}$. It is not hard to show that $\keep_i=1-o(1)$ (see~(\ref{ekk}) below).  So for all relevant values of $i$, $L_i$ is positive and $\D^{2/3}$ is a second-order term in~(\ref{erec1}) and~(\ref{erec2}).}

We will prove:

\begin{lemma}\label{lle} For every $i\geq 0$, every edge $e$ that is uncoloured at the beginning of iteration $i$, \mm{each endpoint $v$ of $e$}, and every $c\in L_i(e)$: if~(\ref{etrunc}) holds for iteration $i$ and $L_i,T_i>\D^{9/10}$ then with probability at least $1-\D^{-10}$, at the beginning of iteration $i+1$  we will have
\begin{enumerate}
\item[(a)] $|L(e)|\geq L_{i+1}$; and
\item[(b)]  $|T(e,v,c)|\leq T_{i+1}$.
\end{enumerate}
\end{lemma}

The Lov\'asz Local Lemma then implies that, with positive probability, the conditions of Lemma~\ref{lle} hold simultaneously for every such $e,c$ and so:

\begin{lemma}\label{llocal}
If~(\ref{etrunc}) holds for iteration $i$ and $L_i,T_i>\D^{9/10}$ then with positive probability~(\ref{etrunc}) holds for iteration $i+1$.
\end{lemma}

\proofstart  
For each edge $e=uv$ and colour $c\in L_i(e)$, we define $A(e)$ to be the event that $|L(e)|<L_{i+1}$ at the beginning of iteration $i+1$, and $B(e,v,c)$ to be the event that $|T(e,v,c)|>T_{i+1}$ at the beginning of iteration $i+1$. If none of these events hold, then~(\ref{etrunc}) holds for iteration $i+1$.
 
Lemma~\ref{lle} says that the probability of each such event is at most $p:=\D^{-10}$.  $A(e)$ is determined by colour assignments and equalizing coin flips for edges incident with $e$; $B(e,v,c)$ is determined by colour assignments and equalizing coin flips for edges within distance two of $e$.  So each event is mutually independent of all events involving edges at distance greater than four, and  thus is mutually independent of all but at most $d:=2\D^4 L_i<\D^5$ other events (see e.g.\ the Mutual Independence Principle in~\cite{mrbook}).   Since \mm{$pd<\inv{4}$} for large $\D$, the Local Lemma completes the proof.
\proofend

A simple analysis of our recursive equations shows that $T_i$ decreases more quickly than $L_i$, and so eventually their ratio will be large enough to allow us to apply Lemma~\ref{ldpcc}. \mm{We must show that this happens before $L_i<\D^{9/10}$, as our procedure stops running if $L_i$ drops below this value.}

\begin{lemma}\label{lrec}
\mm{For every sufficiently small $\e>0$,} there is an $X=X(\e)$ such that for $I=X\ln\D$ we have $L_I>10 T_I$ \mm{and $L_I,T_I>\D^{9/10}$}.
\end{lemma}

\proofstart  \mm{Note that $L_1/T_1=1+\e$. We will prove inductively that $L_i/T_i$ increases with $i$.

Our first useful bound is: If $L_i/T_i\geq 1+\e$ then:
\begin{equation}\label{ekk}
1\geq\keep_i\geq 1-\frac{T_i}{L_i\ln\D}> 1-\frac{1}{(1+\e)\ln\D}.
\end{equation}

Therefore for any constant $X$ and $i\leq X\ln\D$, if $L_i/T_i\geq 1+\e$ and $L_j,T_j\geq \D^{9/10}$ for all $0\leq j\leq i$ then:

 \begin{eqnarray}
\nonumber  L_i>T_i &>& T_1\prod_{j=1}^{i-1}\left(1-\frac{1-\e/2}{\ln\D}\keep^2_j\right)\keep_j\\
& >&T_1\left( 1-\frac{1}{\ln\D}\right)^{I}\left( 1-\frac{1}{(1+\e)\ln\D}\right)^{I}>\D e^{-2X}>\D^{9/10},\label{eti}
 \end{eqnarray}
 for $\D$ sufficiently large in terms of $X,\e$.
 
 We will prove that for any constant $X$, if for all $0\leq j\leq i\leq X\ln\D$ we have $L_j/T_j\geq 1+\e$ and $L_j,T_j\geq \D^{9/10}$  then
\begin{equation}\label{einduct}
\frac{L_{i+1}}{T_{i+1}}\geq \frac{L_i}{T_i}\times\left(1+\frac{\e}{\mm{4}\ln\D}\right),
\end{equation}
for $\D$ sufficiently large in terms of $X,\e$.   It follows inductively that for all $1\leq i\leq I=X\ln \D$ we have $L_i/T_i\geq 1+\e$ and, by~(\ref{eti}), $L_i,T_i\geq \D^{9/10}$.  So the bound in~(\ref{einduct}) holds for all $1\leq i\leq I=X\ln \D$.

To prove~(\ref{einduct}), we first establish bounds on our recursive equations for $L_i,T_i$.  The assumptions that $L_i/T_i>1+\e$ (and so~(\ref{ekk}) holds) and $L_i,T_i>\D^{9/10}$ imply:}
\begin{eqnarray}
L_{i+1}&=&L_i\times \keep^2_i - \D^{2/3}>L_i\times \keep^2_i (1- \D^{-1/5})  \label{eqa}\\
\nonumber T_{i+1}&=& T_i\times\left(1-\frac{\mm{1-\e/2}}{\ln\D}\keep^2_i\right)\times\keep_i + \D^{2/3}\\
&<&T_i\times\left(1-\frac{\mm{1-\e/2}}{\ln\D}\keep^2_i\right)\times\keep_i (1 + \D^{-1/5}). \label{eqb}
\end{eqnarray}
\mm{ Therefore}
\begin{eqnarray*}
\frac{L_{i+1}}{T_{i+1}}&\geq&\frac{L_i}{T_i}\times\frac{\keep_i}{1-\frac{\mm{1-\e/2}}{\ln\D}\keep^2_i}\times\frac{1- \D^{-1/5}}{1+ \D^{-1/5}}\\
&>&\frac{L_i}{T_i}\times\left(1-\frac{1}{(1+\e)\ln\D}\right)\times \left(1+\frac{\mm{1-\e/2}}{\ln\D}\keep^2_i\right)\times\left(1- 2\D^{-1/5}\right)\\
&>&\frac{L_i}{T_i}\times\left(1-(1-\e+\e^2)\frac{1}{\ln\D}\right)
\times \left(1+\frac{\mm{1-2\e/3}}{\ln\D}\right)\qquad\mbox{\mm{by (\ref{ekk})}}\\
&>&\frac{L_i}{T_i}\times\left(1+\frac{\e}{\mm{4}\ln\D}\right),
\end{eqnarray*}
for $\e<\mm{\inv{12}}$ and $\D$ sufficiently large. \mm{This establishes~(\ref{einduct}). Therefore, if} $I=X\ln\D$ where $X$ is a constant that is sufficiently large in terms of $\e$,
\[\frac{L_I}{T_I}>\frac{L_0}{T_0}\times\left(1+\frac{\e}{\mm{4}\ln\D}\right)^I>(1+\e)\times\left(1+\frac{\e}{\mm{4}}\right)^X>10.\]
\mm{This and~(\ref{eti}) prove the lemma.}
\proofend

Our main theorem follows immediately:

\noindent{\bf Proof of Theorem~\ref{mt}:} Setting $I=X\ln\D$ as in Lemma~\ref{lrec}, (\ref{eti}) says that $T_I>\D^{9/10}$ for  $\D$ sufficiently large in terms of $\e$, so our procedure runs for at least $I$ iterations. Since $L_i\geq T_i$, and by the looping rule of our procedure, we have $L_i,T_i>\D^{9/10}$ at every iteration. So Lemma~\ref{llocal} shows inductively that with positive probability~(\ref{etrunc}) holds at the beginning of  every iteration. Thus Lemma~\ref{lrec} and the fact that $L_i/T_i$ is increasing (as shown in the proof of Lemma~\ref{lrec}) yields that with positive probability, when the algorithm terminates, we will have $|L(e)|\geq L_i$ and $|T(e,v,c)|\leq T_i$ for every uncoloured edge $e$, endpoint $v$ of $e$  and colour $c\in L(e)$ where $L_i>10 T_i$.  Now Lemma~\ref{ldpcc} shows that we can complete the colouring. 

This establishes that for every $\e>0$, there exists $\D(\e)$ such that every graph of maximum degree  $\D\geq\D(e)$ has edge correspondence number at most $(1+\e)\D$.  This implies our main theorem.
\proofend

\section{Concentration}\label{sconc}
In this section, we prove our concentration lemma:

\noindent {\bf Proof of  Lemma~\ref{lle}:}\\
{\em Part (a):}  For any colour $c$ remaining in $L(e)$ after step 2(a) of iteration $i$, the probability that $c$ is not removed from $L(e)$ during the remaining steps of iteration $i$ is exactly \mm{$\keep^2_i$}, as explained in section~\ref{srcp}.  

\mm{
{\em Observation 1:} The event that $c$ is not removed from $L(e)$ is mutually independent of the  corresponding events for any other colours of $L(e)$.  

This follows immediately from: (i) for every edge $f\sim e$ and $c'\in L(f)$ there is at most one $c\in L(e)$ such that $c':f$ blocks $c:e$, and (ii) whether $c'$ is assigned to $f$ is independent of the choice to assign any other colour to $f$ or to assign any colour to any other edge.  
}

So the number of colours remaining after those steps is distributed like ${\rm Bin}(L_i,\keep^2_i)$. \mm{Our hypothesis states $L_i>\D^{9/10}$ and we know $\keep_i=1-o(1)$ by~(\ref{ekk}). So the} Chernoff Bounds imply the probability that \mm{fewer} than $\mm{L_{i+1}}=L_i\times\mm{\keep^2_i} - \D^{2/3}$ \mm{colours remain} is at most
\[2e^{-\D^{4/3}/3 L_i\mm{\keep_i^2}}<\D^{-11},\]
for large $\D$ since \mm{$L_i<2\D$ and $\keep_i<1$}.

{\em Part (b):}  

\mmm{
As described earlier, we will show that $|{T'}_{i+1}(e,v,c)|$ is concentrated, and thus is less than $T_{i+1}$ with sufficiently high probability.  This suffices since $|{T}_{i+1}(e,v,c)|\leq |{T'}_{i+1}(e,v,c)|$.

Recall that an edge $f=vw\in T_i(e,v,c)$ is not  in ${T'}_{i+1}(e,v,c)$ if (a) $f$ is assigned and keeps a colour, or (b) $L(f)$ loses $c'$ at $w$, where $c'$ is the unique colour in $L(f)$ such that $c':f$ blocks $c:e$.

We define:  $X_1$ is the number of edges $f=vw\in T_i(e,v,c)$ such that $f$ does not keep a colour.  $X_2$ is the number of edges $f=vw\in T_i(e,v,c)$  such that $f$ does not keep a colour and $L(f)$ loses $c'$ at $w$, where $c'$ is the unique colour in $L(f)$ such that $c':f$ blocks $c:e$.  Thus:

\begin{equation}\label{etxxx}
|T_{i+1}'(e,v,c)|= X_1-X_2.
\end{equation}

We will apply Talagrand's Inequality as stated in Section~\ref{spt} to prove that $X_1, X_2$ are both concentrated. The independent trials will be: (i) for each edge $f$ incident to $v$ or a neighbour of $v$, the set of colours that are assigned to $f$ in Step 2(b.i) and (ii) for each edge $f$ incident to $v$, the set of equalizing coin flips involving $f$.  Note that these trials determine $X_1,X_2$.  

 To be clear: exposing the set of colours assigned to an edge $f$ is one trial, not $|L(f)|$ trials, and exposing the outcomes of all equalizing coinflips involving $f$ is also one trial. So there are a total of at most $\D^2$ trials.

To show that $X_1$ and $X_2$ are concentrated, we will apply Talagrand's Inequality, with
\[D=2\D^{1.2}.\]

Consider any outcome $y$ of our trials, i.e. assignments of colours and equalizing coin flips. 

For every edge $f$ counted by $X_1(y)$, we choose a set of trials $T(f)$ which certify that $f$ is counted by $X_1(y)$; if $f$ is also counted by $X_2(y)$ then $T(f)$ will also certify this. We place the colour assignments to $f$ in $T(f)$; we place the equalizing coinflips for $f$ in $T(f)$; and for each colour $\alpha$ assigned to $f$ that was removed from $f$ because of the colour assignments to another edge,  we chose one such edge $f'$ and place the colour assignments for $f'$ into $T(f)$.   
If $f=vw$ is also counted by $X_2(y)$ then we place into $T(f)$ the colour assignments to an edge incident with $w$ that caused $L(f)$ to lose the colour blocking $c:e$. (Perhaps that trial was already in $T(f)$ if that colour was assigned to $f$.)

 For each trial $T_j$,  if $T_j$ is a set of equalizing coinflips then we set $b_j=1$.   If $T_j$ is the assignment of colours to an edge $f$ incident with $v$ then we set $b_j=\D^{1/10}+1$.   If $T_j$ is the assignment of colours to an edge  that is not incident to $v$ but that is in $T(f)$ for some $f$ counted by $X_2$, then we set $b_j=2$.

We need to bound the probability of the event $F$ from our statement of Talagrand's Inequality.  

We define the event $Q$ to be the event that:\\
\mbox{  }\\
{\em \noindent (Q1) every edge $f\in T_i(e,v,c)$ is assigned fewer than $\D^{1/10}$ colours; and\\ 
\noindent (Q2) for every edge $f'$ incident to $v$, there are fewer than $\D^{1/10}$ edges  $f\in T_i(e,v,c)$ such that $f,f'$ receive colours that block each other; i.e. $f$ receives a colour $\a$ and $f'$ receives a colour $\a'$ such that $\a:f$ blocks $\a':f'$
} 

We will show that $Q$ implies $F$.

First we bound $\sum b_i^2$.
There are at most $\D$ uncoloured edges incident to $v$, and so at most $\D$ trials $T_j$ with $b_j=\D^{1/10}+1$ and at most $\D$ trials $T_j$ with $b_j=1$.   If $(Q1)$ holds, then each edge $f$ counted by $X_1(y)$ or $X_2(y)$ has  fewer than $\D^{1/10}$ assigned colours and so the colour assignments to at most $\D^{1/10}$ edges  incident to $f$ but not $v$ are in $T(f)$; i.e. one edge for each of the colours assigned to $f$ and possibly one edge for the colour  in $L(f)$ that blocks $c:e$. This yields at most $\D\times\D^{1/10}$ trials $T_j$ with $b_j=2$. Thus for sufficiently large $\D$,
\[\sum b_i^2<\D\times (\D^{1/10}+1)^2 + \D\times 1^2 + \D\times\D^{1/10}\times 2^2<   2\D^{1.2}=D,\]
 as required for $F$.

Let $z$ be any other possible outcome of the trials. 
If $f$ is counted by $X_1(y)$ but not by $X_1(z)$, or if $f$ is counted by $X_2(y)$ but not by $X_2(z)$ then $y,z$ must differ on some trial in $T(f)$.     If $(Q2)$ holds then each trial consisting of colour assignments to an edge $f'$ incident to $v$ can only be in $T(f)$ for at most $\D^{1/10}+1$ choices of $f$ (including $f=f'$).  Because $G$ is simple, each trial consisting of the colour assignments to an edge not incident with $v$ can only be in $T(f)$ for at most 2 choices of $f$; the at most two edges from $v$ to an endpoint of that edge. A trial consisting of the equalizing coin flips for an edge can only be in $T(f)$ for that edge $f$. It follows that
\begin{eqnarray*}
X_1(z)&\geq& X_1(y)-\sum_{y_j\neq z_j} b_j\\
X_2(z)&\geq& X_2(y)-\sum_{y_j\neq z_j} b_j
\end{eqnarray*}
as required for $F$.

So $Q$ implies $F$ and thus $\pr(\ov{F})\leq \pr(\ov{Q})$ for both $X_1$ and $X_2$.

Straightforward calculations show that $\pr(\ov{Q})$ is very small.  Indeed, the expected number of edges in $T(e,v,c)$ which receive more than $\D^{1/10}$ colours is at most 
\[T_i{L_i\choose\D^{1/10}}\left(\inv{L_i\ln\D}\right)^{\D^{1/10}}
<T_i\left(\frac{eL_i}{\D^{1/10}}\right)^{\D^{1/10}}\left(\inv{L_i\ln\D}\right)^{\D^{1/10}}
=T_i\left(\frac{e}{\D^{1/10}\ln\D}\right)^{\D^{1/10}}<\hf\D^{-11}.\]

For the expected number of edges $f'$ violating (Q2): First we choose $f'$ incident with $v$.  We must pick  $\D^{1/10}$ edges $f\in T(e,v,c)$, and for each such $f$ a colour assigned to $f$ that blocks a colour assigned to $f'$.   We start by choosing the colours assigned to $f'$.  Note that each such colour could block colours assigned to multiple choices for $f$. We choose integers $\ell, a_1,...,a_{\ell}>0$ where $\sum_{j=1}^{\ell}a_j=\D^{1/10}$, and then for each $j=1,...,\ell$, we select a colour $\a_j\in L(f')$ which will block colours assigned to $a_j$ edges $f$.  We then select those $a_j$ edges from the at most $T_i$ edges in $T(f',v,\a_j)\cap T(e,v,c)$; note that for each selected edge $f$, there is exactly one choice for the colour assigned to $f$ which blocks $f':\a_j$, by the definition of correspondence colouring.  We multiply by the probability that each of these $\D^{1/10}+\ell$ colour-assignments occur.

 Putting it all together, the expected number of violations to (Q2) is at most: 
\begin{eqnarray*}
&&\D\times\sum_{\ell\geq 0}\sum_{a_1,...,a_{\ell}>0; \sum a_j=\D^{1/10}}
L_i^{\ell}\left(\prod_{j=1}^{\ell}{T_i\choose a_j}\right)\left(\inv{L_i\ln\D}\right)^{\D^{1/10}+\ell}\\
&<&\D\times\sum_{\ell\geq 0} \sum_{a_1,...,a_{\ell}>0; \sum a_j=\D^{1/10}}L_i^{\ell}\left(\prod_{j=1}^{\ell}T_i^{a_j}\right)\left(\inv{L_i\ln\D}\right)^{\D^{1/10}+\ell}\\
&<&\D\times\sum_{\ell\geq 0} 2^{\D^{1/10}}L_i^{\ell}T_i^{\D^{1/10}}\left(\inv{L_i\ln\D}\right)^{\D^{1/10}+\ell}\\
&&\qquad\qquad\qquad\mbox{as the number of choices for $\ell,a_1,...,a_{\ell}$ is at most $2^{\D^{1/10}}$}\\
&=&\D\left(\frac{T_i}{L_i\ln\D}\right)^{\D^{1/10}}\sum_{\ell\geq 0}\left(\inv{\ln\D}\right)^{\ell}\\
&<&\hf\D^{-11}
\end{eqnarray*}
for $\D$ sufficiently large, since we showed in the proof of Lemma~\ref{lrec} that $L_i/T_i>1+\e>1$.
So by Markov's Inequality,
\begin{equation}\label{eprf}
\pr(\ov{F})\leq\pr(\ov{Q})<\D^{-11}.
\end{equation}

We now apply  Talagrand's Inequality with $D=2\D^{1.2}, \xmax=\D$ and $t=\inv{4} \D^{2/3}$.  First note that
$t+35\sqrt{\D}+\xmax\pr(\ov{F})<\hf \D^{2/3}$.  So for $i=1,2$:
\begin{eqnarray*}
\pr(|X_i - \ex(X_i)|>\hf\D^{2/3}) &\leq& \pr(|X_i - \ex(X_i)|>t+35\sqrt{\D}+\xmax\pr(\ov{F}))\\
&<& 2\pr(\ov{F})+5e^{-\frac{t^2}{4D}} < \hf\D^{-10}.
\end{eqnarray*}
Recalling~(\ref{et'}),~(\ref{erec2}), and~(\ref{etxxx}), if $|T^{'}_{i+1}(e,v,c)|>T_{i+1}$ then at least one of $|X_1 - \ex(X_1)|, |X_2 - \ex(X_2)|$ must be greater than $\hf\D^{2/3}$.  This proves part (b) for $\D$ sufficiently large.
\proofend

}

\section{Hypergraphs}

We close by remarking that our main theorem also holds for linear $k$-uniform hypergraphs for $k=O(1)$.  I.e., for any constant $k$, and any hypergraph $H$ where every hyperedge contains exactly $k$ vertices, every pair of vertices lies in at most one hyperedge,  and every vertex lies in at most $\D$ hyperedges, we have $\chi'_{DP}(G)=\D+o(\D)$.  The proof is a very straightforward adaptation of the proof of Theorem~\ref{mt}. \mm{We outline it here: 

We first remark that our definintion of $\chi'_{DP}$ extends naturally to linear hypergraphs.  

In the statement of Lemma~\ref{ldpcc}, the constant 8 changes to $4k$.  So in the halting condition of our procedure, $L_i>10T_i$ becomes $L_i>5k T_i$, and the constant 10 is replaced with $5k$ as appropriate throughout the proof.

For each hyperedge $e$, we  define $T_i(e,v,c),T'_i(e,v,c),\eq_i(e,v,c)$ for each of the $k$ vertices $v$ in $e$. 

Equalizing coinflips ensure that for every $v\in e$ and $c \in L(e)$, the probability that $L(e)$ loses $c$ at $v$ is exactly $1-\keep_i$. So now the probability that $c$ remains in $L(e)$ is $\keep_i^k$.  So our recursive equation for $L_i$ becomes:
\[L_{i+1}= L_i\times\keep_i^{k}-\D^{2/3}.\]

This time, $T'_{i+1}(e,v,c)$ is the number of edges $f\subseteq T_{i}(e,v,c)$ such that during iteration $i$, (a) $f$ does not retain a colour  and (b) $L(f)$ does not lose $c'$ at {\em any of its $k-1$ vertices other than $v$}, where $c'$ is the unique colour in $L(f)$ such that $c':f$ blocks $c:e$. So equation~(\ref{et'}) becomes:
\[\ex[ |T'_{i+1}(e,v,c)|]< |T_{i}(e,v,c)| \times\left(1-\frac{1-\e/2}{\ln\D}\keep^k_i\right)\times\keep^{k-1}_i.
\]
So our recursive equation for $T_i$ becomes:
\[T_{i+1}= T_i\times\left(1-\frac{1-\e/2}{\ln\D}\keep^k_i\right)\times\keep^{k-1}_i+\D^{2/3}.\]

The statement of Observation~\ref{oZY} is modified to include $k$ $Y$-events. The fact that $H$ is linear ensures that this Observation still holds.

The proof of Lemma~\ref{lle} changes only in very straightforward places,\mmm{ for example: For each hyperedge $f'$ not incident to $v$, such that $f'$ is in some $T(f)$, we set $b_j=k$ for the trial $T_j$ consisting of the colour assignments to  $f'$. This is because $T_j$ can be in $T(f)$ for up to $k$ hyperedges $f$ linking $v$ to $f'$; the fact that $H$ is linear is important here.}  The same is true for Lemma~\ref{llocal}. The change in our recursive equations results in very minor changes to the calculations in the proof of Lemma~\ref{lrec}.  Those three lemma statements remain the same, except for changing 10 to $5k$ in Lemma~\ref{lrec}.

This modified proof yields:

\begin{theorem}\label{mt2} For any constant $k$, let $H$ be any linear $k$-uniform hypergraph with maximum degree $\D$.  Then $\chi'_{DP}(H)=\D+o(\D)$.
\end{theorem}

}

 \section*{Acknowledgement}  We are grateful to Runrun Liu, Stijn Cambie, and an anonymous referee for providing several corrections and improvements to earlier drafts.

\section*{\mmm{Appendix: Talagrand's Inequality}}

Here, we show how to obtain our statement of Talagrand's Inequality from Talagrand's original statement.  We remark that our statement is very similar to ``Talagrand's Inequality V'' from Chapter 20 of~\cite{mr3}. Unfortunately, there is an error in that statement (specifically, there is an error in Fact 20.1 of that book) and so we prove this version instead.  That error also discussed in~\cite{kp}.   The statement we use here is derived from the statement used in~\cite{mr3} and the proof is almost identical.

We start with Talagrand's original statement.   Consider $n$ independent random trials $T_1,...,T_n$, and let $\cala$ be the set of all possible sequences of outcomes.  Consider any subset $A\subseteq\cala$ and any real $\ell\geq 0$.   We define $A_{\ell}$ as follows: A set of outcomes $y=(y_1,...,y_n)\in\cala$ is defined to be in $A_{\ell}$ if for every set of real $b_1,...,b_n$ there is at least one $z=(z_1,...,z_n)\in A$ such that
\begin{equation}\label{etsum}
\sum_{y_i\neq z_i} b_i < \ell\left(\sum_{i=1}^n b_i^2\right)^{1/2}.
\end{equation}
 Talagrand's original statement (see Theorem 4.1.1 of \cite{mt}) is:
\begin{theorem} For any $n$ independent trials $T_1,...,T_n$, any $A\subseteq \cala$ and any real $\ell$,
\[\pr(A)\times\pr(\ov{A_{\ell}}) \leq \rme^{-\ell^2/4}.\]
\end{theorem}
We can think of $A_{\ell}$ as the outcomes that are in some sense close to $A$, i.e.  close to at least one point in $A$. So for sets $A$ that are reasonably large; i.e. $\pr(A)$ is reasonably large, the probability of our outcome being far from $A$ is exponentially small.

Very roughly speaking, we will set things up so that if $X$ is far from $\Exp(X)$ then our outcomes must be far from a reasonably large set $A$, thus showing that the probability of $X$ being far from $\Exp(X)$ is exponentially small.  In fact, we will first work with the median, $\med(X)$, instead so we can use the fact  that the probability of $X$ being at most $\med(X)$ is reasonable large.  After showing that $X$ is concentrated around $\med(X)$ then we will complete the proof by showing that $\Exp(X)$ is close to $\med(X)$; this is where the extra  $35\sqrt{D}+\xmax\times\pr(\ov{F})$  term below comes from.

Before getting into the details, we recall the statement from section~\ref{spt}:

{\em Let $X$ be a random variable determined by the
independent trials $T_1,...,T_n$. Let $\xmax,D>0$ and suppose that $X$ always satisfies $0\leq X\leq\xmax$. Let $F$ be the event that for the outcome $y=(y_1,...,y_n)$ of the trials there exists $b_1,....,b_n>0$ such that
\begin{enumerate}
\item[(i)] $\sum_{i=1}^n b_i^2\leq D$; and
\item[(ii)] for any possible outcome $z=(z_1,...,z_n)$ of the trials, we have
\[X(z)\geq X(y)-\sum_{y_i\neq z_i} b_i.\]
\end{enumerate}
Then for any $t\geq 0$ we have
\begin{equation}\label{tigoal}
\pr\left(|X-\ex(X)|>t +  35\sqrt{D}+\xmax\times\pr(\ov{F})\right)\leq 2\pr(\ov{F}) + 5e^{-\frac{t^2}{4D}}.
\end{equation}
}

\proofstart
Note that we can assume
\begin{equation}\label{ef12}
\pr(\ov{F})<\hf,  e^{-\frac{t^2}{4D}}<\inv{5},
\end{equation}
else (\ref{tigoal}) is trivial.
We first bound the probability of $X$ being much higher than $\med(X)$.  We define $A=\{x:X(x)\leq\Med(X)\}$, 
and we set $\ell=t/\sqrt{D}$. Thus $\pr(A)\geq\hf$, and thus Talagrand's Inequality implies $\pr(\ov{A_{\ell}})$ is small.  

We define $C=\{y:X(y)>\Med(X)+t\}\cap F$. \footnote{Adding ``$\cap F$'' to the definition of $C$ and then keeping track of the minor implications is the only significant difference between this proof and the corresponding proof in~\cite{mr3}.} To prove that $\pr(C)$ is small, we will show
\begin{equation}\label{ecal1}
C\subseteq\ov{A_{\ell}}.
\end{equation}
Consider any $y=(y_1,...,y_n)\in C$.  We wish to prove $y\notin A_{\ell}$, i.e. that there is some $b_1,....,b_n$ such that~(\ref{etsum}) fails to hold for all $z\in A$.  Since $y\in C\subseteq F$ we can take the values $b_1,...,b_n$ from the definition of $F$.  
 
Consider any $z\in A$. Note that, by the definitions of $A,C$, we have $X(y)-X(z)>t$.  Therefore, condition (ii) of the definition of $F$ yields
\[ \sum_{y_i\neq z_i} b_i\geq X(y)-X(z)>t.\]
Condition (i) says $\sum_{i=1}^n b_i^2\leq D=(t/\ell)^2$.  Putting these together shows that~(\ref{etsum}) does not hold; i.e. $\sum_{y_i\neq z_i} b_i > \ell\left(\sum_{i=1}^n b_i^2\right)^{1/2}$.  Since this is true for every $z\in A$, this shows $y\notin A_{\ell}$, thus proving~(\ref{ecal1}). Sine $\pr(A)\geq\hf$, Talagrand's Inequality implies
\begin{equation}\label{ec1}
\pr(C)\leq\pr(\ov{A_{\ell}})\leq2e^{-\ell^2/4}=2e^{-t^2/4D}.
\end{equation}
To bound the probability that $X$ is much smaller than $\med(X)$, we set $A'=\{x:X(x)\geq\Med(X)\}\cap F$ and $C'=\{y:X(y)<\Med(X)-t\}$.  We repeat the same argument as above, but with $A'$ playing the role of $C$ and vice versa since now $A'$ contains the outcomes with the higher values of $X$.  We omit the repetitive details which yield:  
\[A'\subseteq\ov{C'_{\ell}}.\]
Therefore, $\pr(\ov{C'_{\ell}})\geq\pr(A')\geq\hf-\pr(\ov{F}))$ and so Talagrand's Inequality and~(\ref{ef12}) gives
\begin{equation}\label{ec2}
\pr(C')\leq e^{-\ell^2/4}/(\hf-\pr(\ov{F}) <e^{-t^2/4D} (3+3\pr(\ov{F})).
\end{equation}
Putting this together and again applying~(\ref{ef12}), we have
\begin{equation}\label{emedconc}
\pr(|X-\med(X)|>t )\leq \pr(C)+\pr(\ov{F})+\pr(C')< 2\pr(\ov{F}) + 5e^{-\frac{t^2}{4D}}.
\end{equation}

To show that this implies concentration around the {\em mean},
we prove that the mean and median do not differ by very much:

\begin{claim}\label{tal2.f1}
Under the preconditions of (\ref{tigoal}),
$|\ex(X)-\med(X)|\leq 35\sqrt{D}+ \xmax\times\pr(F)$.
\end{claim}

This will prove our statement of Talagrand's Inequality as if $|X-\ex(X)|>t +  35\sqrt{D}+\xmax\times\pr(\ov{F})$ then Claim~\ref{tal2.f1} implies that $|X-\med(X)|>t $ and so~(\ref{tigoal}) follows from~(\ref{emedconc}).  

{\bf Proof of Claim~\ref{tal2.f1}.}  Note that 
\[|\ex(X)-\med(X)|\leq |\ex(X|F)-\med(X)|\times\pr(F) + X_{max}\times\pr(\ov{F}).\]
So most of our work will be to bound $|\ex(X|F)-\med(X)|$.
First, observe that $\ex(X|F)-\med(X)=\ex(X-\med(X)|F)$.
We will bound the absolute value of this latter term by partitioning the positive
real line into the intervals
$I_i=(i\times \sqrt{D},(i+1)\times \sqrt{D}]$,
defined for each integer $i\geq0$. Clearly, $\left|\ex(X-\med(X)|F)\right|\times\pr(F)$,  is at most the sum over all $I_i$ of
the maximum value in $I_i$ times the probability of the event $(|X-\med(X)|\in I_i)\cap F$, which is
\begin{eqnarray*}
&&\sum_{i\geq0}(i+1)\times \sqrt{D}\times\pr((|X-\med(X)|\in I_{i}) \cap F)\\
&=&\sum_{i\geq0}\sqrt{D}\times\pr((|X-\med(X)|\in \cup_{j\geq i}I_j)\cap F).
\end{eqnarray*}
Now recall that $C=\{y:X(y)>\Med(X)+t\}\cap F$ and $C'=\{y:X(y)<\Med(X)-t\}$. So
setting $t=i\sqrt{D}$ and applying~(\ref{ec1}),~(\ref{ec2}) and~(\ref{ef12}) yields
 a bound on  $\pr((|X-\med(X)|\in \cup_{j\geq i}I_j)\cap F)$ of:
\[(5+3\pr(\ov{F}))\rme^{-{i^2/4}}<7\rme^{-{i^2/4}}.\]
Therefore:
\[|\ex(X-\med(X)|F)|\times\pr(F)
<7\sqrt{D}\times\sum_{i\geq0}\rme^{-{i^2/4}}.\]
It is straightforward to bound 
$\sum_{i\geq0}\rme^{-{i^2\over 4}}< 4+\sum_{i\geq4}\rme^{-i}<5$.  Therefore
\[|\ex(X)-\med(X)|< 35\sqrt{D} + X_{max}\times\pr(\ov{F}),\]
as required
\proofend

\end{document}